\documentclass[11pt]{article}
\usepackage{geometry}                
\usepackage{amsmath,amssymb}
\usepackage{amsthm}
\usepackage{hyperref}
\usepackage{graphicx,epsfig,color}
\usepackage{booktabs,bm,multirow}
\usepackage{cases}

\newtheorem{theorem}{Theorem}

\newtheorem{lemma}[theorem]{Lemma}

\input{siam.sty}

\title{Stochastic gradient descent for linear least squares problems with partially observed data}
\author{Kui Du\thanks{School of Mathematical Sciences, Xiamen University, Xiamen 361005, China ({\tt kuidu@xmu.edu.cn}).},\quad Xiao-Hui Sun\thanks{School of Mathematical Sciences, Xiamen University, Xiamen 361005, China ({\tt 19020190154621@stu.xmu.edu.cn}).}} 
\date{}

\begin{document}
\maketitle

\begin{abstract}  
\vspace{.5mm} 

We propose a novel stochastic gradient descent method for solving linear least squares problems with partially observed data. Our method uses submatrices indexed by a randomly selected pair of row and column index sets to update the iterate at each step. Theoretical convergence guarantees in the mean square sense are provided. Numerical experiments are reported to demonstrate the theoretical findings.

{\bf Keywords}. linear least squares problem, partially observed data, stochastic gradient descent

{\bf AMS subject classifications}: 65F10, 65F20, 68W20\end{abstract}

\section{Introduction} 
In recent years, stochastic iterative methods for solving large-scale linear equations or linear least squares problems have been greatly developed due to low memory footprints, such as randomized Kaczmarz \cite{strohmer2009rando}, randomized coordinate descent \cite{leventhal2010rando}, and their extensions, e.g., \cite{zouzias2013rando,ma2015conve,gower2015rando,richtarik2017stoch,bai2018greed,necoara2019faste,bai2019greed,zhang2019new,liu2019varia,wu2020proje,du2019doubl,du2020pseud,chen2020error}.
However, the performance of these methods are often limited \cite{ma2019stoch} when solving the problems with partially observed data, which can arise due to lost of data, errors in data recording, or cost of data acquisition, etc.

In this paper we consider solving the linear least squares problem \begin{equation}\label{ls}\min_{\mbf x\in\mbbr^n}\|{\bf Ax-b}\|_2,\end{equation} where $\mbf A\in\mbbr^{m\times n}$ ($m\geq n$) has full column rank and $\|\cdot\|_2$ denotes the Euclidean norm. The least squares solution of this problem is $\mbf A^\dag\mbf b$, where $\mbf A^\dag$ is the Moore-Penrose generalized inverse \cite{ben2003gener}. Sometimes the matrix $\mbf A$ and the vector $\mbf b$ are partially observed, i.e., only partial entries of $\mbf A$ and $\mbf b$ are available. As a concrete example, suppose $\mbf A$ is the score matrix of $m$ users for $n$ services, and $\mbf b$ contains the $m$ weighted total scores from each user for these services. Each user may not give scores for all the $n$ services, or may not give a weighted total score for these $n$ services, but the survey company wants to know the weight of each service in the weighted total score. That is to say, we need to solve the linear least squares problem (\ref{ls}) with only partial entries of $\mbf A$ and $\mbf b$ are available.

Let $\{\delta_{ij}\}_{i=1,j=1}^{m,n}$ denote independent and identically distributed (i.i.d.) Bernoulli random variables satisfying $$\delta_{ij}=\begin{cases}1, & \mbox{with  probability}\ p,\\ 0, &\mbox {with  probability}\ 1-p,\end{cases}$$ and $\{\delta_{i}\}_{i=1}^{m}$ denote i.i.d. Bernoulli random variables satisfying $$\delta_{i}=\begin{cases}1, & \mbox{with  probability}\ q,\\ 0, &\mbox {with  probability}\ 1-q.\end{cases}$$ If we use $\delta_{ij}=1$ or $\delta_i=1$ to indicate the availability of an element in $\mbf A$ or $\mbf b$, and $\delta_{ij}=0$ or $\delta_i=0$ to indicate a missing entry, then the partially observed data are \begin{align}\label{pd} \wh{\mbf A}=[\delta_{ij}\mbf A_{ij}]_{i=1,j=1}^{m,n}\quad  \mbox{ and }\quad  \wh {\mbf b}=[\delta_i\mbf b_i]_{i=1}^m.\end{align} So the linear least squares problem with partially observed data is: \beq\label{pls}\mbox{Given}\quad \wh{\mbf A}, \wh{\mbf b},\quad \mbox{find the unique least squares solution}\quad \mbf A^\dag\mbf b=\argmin_{\mbf x\in\mbbr^n}\|{\bf Ax-b}\|_2.\eeq 

Solving the problems with partially observed data by standard methods, such as data imputation \cite{bradley1994missi}, matrix completion \cite{cai2010singu,keshavan2010matri2,keshavan2010matri,recht2011simple}, and maximum likelihood estimation  \cite{dempster1977maxim,little2002stati}, can be wasteful, produces biases, or is impractical for large-scale problems. Recently, Ma and Needell \cite{ma2019stoch} proposed a stochastic gradient descent (SGD) method for the linear least squares problem (\ref{ls}) with partially observed ${\mbf A}$ as given in (\ref{pd}) and fully observed $\mbf b$. Their method uses a randomly selected row of $\wh{\mbf A}$ to construct a stochastic gradient at each step. In this paper, we consider the more general case, i.e., both $\mbf A$ and $\mbf b$ are partially observed. 
 
{\it Main contributions}. We propose a novel stochastic gradient descent method for solving the linear least squares problem (\ref{ls}) with partially observed data (\ref{pd}) and prove its convergence theory. At each step, the new method uses submatrices indexed by a randomly selected pair of row and column index sets to design a stochastic gradient. Numerical examples are reported to illustrate the convergence of the new method.

{\it Organization of this paper}. In the rest of this section, we give some notation. In Section 2 we describe the proposed stochastic gradient descent method and prove its convergence theory. In Section 3 we report the numerical results. Finally, we present brief concluding remarks in Section 4.

{\it Notation}. For any random variables $\bm\xi$ and $\bm\zeta$, we use $\mbbe\bem\bm\xi\eem$ and $\mbbe\bem\bm\xi\ |\bm\zeta\eem$ to denote the expectation of $\bm\xi$ and the conditional expectation of $\bm\xi$ given $\bm\zeta$, respectively. For an integer $m\geq 1$, let $[m]:=\{1,2,3,\ldots,m\}$.  For any vector $\mbf b\in\mbbr^m$, we use $\mbf b_i$, $\bf b^\rmt$ and $\|\mbf b\|_2$ to denote, the $i$th entry, the transpose and the Euclidean norm of $\mbf b$, respectively. We use $\mbf I$ to denote the identity matrix whose order is clear from the context. For any matrix $\mbf A\in\mbbr^{m\times n}$, we use $\mbf A_{i,j}$, $\mbf A_{i,:}$, $\mbf A_{:,j}$ $\mbf A^\rmt$, $\mbf A^\dag$, $\|\mbf A\|_\rmf$, $\ran(\mbf A)$, and $\sigma_{\rm min}(\mbf A)$ to denote the $(i,j)$ entry, the $i$th row, the $j$th column, the transpose, the Moore-Penrose pseudoinverse, the Frobenius norm, the column space, and the smallest nonzero singular values of $\mbf A$, respectively.   For index sets $\mcali\subseteq[m]$ and $\mcalj\subseteq[n]$, let $\mbf A_{\mcali,:}$, $\mbf A_{:,\mcalj}$, and $\mbf A_{\mcali,\mcalj}$ denote the row submatrix indexed by $\mcali$, the column submatrix indexed by $\mcalj$, and the submatrix that lies in the rows indexed by $\mcali$ and the columns indexed by $\mcalj$, respectively. Similarly, we use $\mbf b_\mcali$ to  denote the column vector consisting of elements of $\mbf b$ indexed by $\mcali$. Given a square matrix $\mbf M$, we denote a matrix containing only the diagonal of $\mbf M$ as $\diag(\mbf M)$. We use $\mbf B\preceq\mbf A$ to denote that $\mbf A-\mbf B$ is positive semidefinite. 

\section{Stochastic gradient descent for partially observed data}

In \cite{du2019doubl} we proposed a doubly stochastic block Gauss-Seidel (DSBGS) algorithm for solving a consistent linear system $\bf Ax=b$. By varying the row partition parameter and the column partition parameter of DSBGS, we recover the randomized Kaczmarz algorithm \cite{strohmer2009rando}, the randomized coordinate descent algorithm \cite{leventhal2010rando}, and the doubly stochastic Gauss-Seidel algorithm \cite{razaviyayn2019linea}. Next we provide a slightly different variant of DSBGS, which will be used to derive our stochastic gradient descent method for solving the problem (\ref{pls}).
 
Let $\{\mcali_1,\mcali_2,\ldots,\mcali_s\}$ denote a partition of $[m]$ such that, for $i,j=1,2,\ldots,s$ and $i\neq j,$ $$ \mcali_i\neq\emptyset,\quad\mcali_i\cap\mcali_j=\emptyset,\quad \bigcup_{i=1}^s\mcali_i=[m].$$ Let $\{\mcalj_1,\mcalj_2,\ldots,\mcalj_t\}$ denote a partition of $[n]$ such that, for $i,j=1,2,\ldots,t$ and $i\neq j,$ $$ \mcalj_i\neq\emptyset,\quad\mcalj_i\cap\mcalj_j=\emptyset,\quad \bigcup_{i=1}^t\mcalj_i=[n].$$ Let $$\mcalp=\{\mcali_1,\mcali_2,\ldots,\mcali_s\}\times \{\mcalj_1,\mcalj_2,\ldots,\mcalj_t\}.$$ To solve the linear least squares problem (\ref{ls}), one approach is to minimize the function 
\beqs\label{fx} f(\mbf x):=\frac{1}{2st}\|\mbf A\mbf x-\mbf b\|_2^2.\eeqs
If a pair of index sets $(\mcali,\mcalj)$ is randomly selected with probability $\frac{1}{st}$, then we obtain \begin{align*} &\hspace{5.2mm}\mbbe\bem \mbf I_{:,\mcalj_j}(\mbf A_{\mcali_i,\mcalj_j})^\rmt(\mbf A_{\mcali_i,:}\mbf x-\mbf b_{\mcali_i})\eem\\ & =\frac{1}{st}\sum_{i=1}^s\sum_{j=1}^t\mbf I_{:,\mcalj_j}(\mbf A_{\mcali_i,\mcalj_j})^\rmt(\mbf A_{\mcali_i,:}\mbf x-\mbf b_{\mcali_i})\\& =\frac{1}{st}\sum_{i=1}^s\sum_{j=1}^t\mbf I_{:,\mcalj_j}(\mbf I_{:,\mcalj_j})^\rmt\mbf A^\rmt\mbf I_{:,\mcali_i}(\mbf I_{:,\mcali_i})^\rmt(\mbf A\mbf x-\mbf b)\\ &=\frac{1}{st}\mbf A^\rmt(\mbf A\mbf x-\mbf b)\\ &= \nabla f(\mbf x).\end{align*} This yields a stochastic gradient descent method  (see Algorithm 1)  for the linear least squares problem (\ref{ls}). Note that Algorithm 1 is a slightly different variant of DSBGS \cite{du2019doubl}, which randomly selects the pair $(\mcali,\mcalj)$ with probability $\|\mbf A_{\mcali,\mcalj}\|^2_\rmf/\|\mbf A\|_\rmf^2.$ 

\begin{center}
\begin{tabular*}{145mm}{l}
\toprule {\bf Algorithm 1:} SGD for the linear least squares problem (\ref{ls})  \\ 
\noalign{\smallskip}\hline \noalign{\smallskip}
\quad Initialize $\mbf x^0\in\mbbr^n$\\ 
\quad {\bf for} $k=1,2,\ldots,$ {\bf do}\\
\quad \hspace{5.2mm} Pick $(\mcali,\mcalj)\in\mcalp$ with probability $\dsp\frac{1}{st}$\\
\quad \hspace{5.2mm} Pick $\alpha_k>0$\\\noalign{\smallskip}
\quad \hspace{5.2mm}  Set $\dsp\mbf x^k=\mbf x^{k-1}-\alpha_k\mbf I_{:,\mcalj}(\mbf A_{\mcali,\mcalj})^\rmt(\mbf A_{\mcali,:}\mbf x^{k-1}-\mbf b_\mcali)$\\
\bottomrule
\end{tabular*}
\end{center}

Directly applying Algorithm 1 to the partially observed data (\ref{pd}), we obtain the update: \beq\label{sgd1}\mbf x^k=\mbf x^{k-1}-\alpha_k\mbf I_{:,\mcalj}(\wh{\mbf A}_{\mcali,\mcalj})^\rmt(\wh{\mbf A}_{\mcali,:}\mbf x^{k-1}-\wh{\mbf b}_\mcali).\eeq 
Note that (see detailed calculations in the proof of Lemma \ref{lem-g}) \begin{align*} &\hspace{5.2mm}\mbbe\bem \mbf I_{:,\mcalj}(\wh{\mbf A}_{\mcali,\mcalj})^\rmt(\wh{\mbf A}_{\mcali,:}\mbf x^{k-1}-\wh{\mbf b}_\mcali)\ | \mbf x^{k-1}\eem \\ & = \frac{p^2}{st} \mbf A^\rmt\mbf A\mbf x^{k-1}+\frac{p-p^2}{st}\diag(\mbf A^\rmt\mbf A)\mbf x^{k-1}-\frac{pq}{st}\mbf A^\rmt\mbf b\\ \nn &\neq \nabla f(\mbf x^{k-1}).\end{align*}  As a result, the iteration (\ref{sgd1}) does not move in the negative gradient direction. Instead of using $\mbf I_{:,\mcalj}(\wh{\mbf A}_{\mcali,\mcalj})^\rmt(\wh{\mbf A}_{\mcali,:}\mbf x^{k-1}-\wh{\mbf b}_\mcali)$, if we use 
\beqs\mbf g(\mbf x^{k-1})=\mbf I_{:,\mcalj}(\wh{\mbf A}_{\mcali,\mcalj})^\rmt\l(\dsp\frac{\wh{\mbf A}_{\mcali,:}\mbf x^{k-1}}{p^2}-\frac{\wh{\mbf b}_\mcali}{pq}\r)-\dsp\frac{1-p}{p^2}\diag\l(\mbf I_{:,\mcalj}(\wh{\mbf A}_{\mcali,\mcalj})^\rmt\wh{\mbf A}_{\mcali,:}\r)\mbf x^{k-1},\eeqs then we have (see Lemma \ref{lem-g})\beq\label{egx}\mbbe\bem \mbf g(\mbf x^{k-1})\ |\mbf x^{k-1}\eem=\frac{1}{st}\mbf A^\rmt(\mbf A\mbf x^{k-1}-\mbf b)=\nabla f(\mbf x^{k-1}).\eeq
This yields a stochastic gradient descent method (see Algorithm 2) for solving the linear least squares problem (\ref{ls}) with partially observed data (\ref{pd}).

\begin{center}
\begin{tabular*}{155mm}{l}
\toprule {\bf Algorithm 2:} SGD for the linear least squares problem with partially observed data (\ref{pd})  \\ 
\noalign{\smallskip}\hline \noalign{\smallskip}
\quad Initialize $\mbf x^0\in\mbbr^n$\\ 
\quad {\bf for} $k=1,2,\ldots,$ {\bf do}\\
\quad \hspace{5.2mm} Pick $(\mcali,\mcalj)\in\mcalp$ with probability $\dsp\frac{1}{st}$\\\noalign{\smallskip}
\quad \hspace{5.2mm} Set $\mbf g(\mbf x^{k-1})=\mbf I_{:,\mcalj}(\wh{\mbf A}_{\mcali,\mcalj})^\rmt\l(\dsp\frac{\wh{\mbf A}_{\mcali,:}\mbf x^{k-1}}{p^2}-\frac{\wh{\mbf b}_\mcali}{pq}\r)-\dsp\frac{1-p}{p^2}\diag\l(\mbf I_{:,\mcalj}(\wh{\mbf A}_{\mcali,\mcalj})^\rmt\wh{\mbf A}_{\mcali,:}\r)\mbf x^{k-1}$ \\ 
\quad \hspace{5.2mm} Pick $\alpha_k>0$\\\noalign{\smallskip}\noalign{\smallskip}
\quad \hspace{5.2mm}  Set $\dsp\mbf x^k=\mbf x^{k-1}-\alpha_k\mbf g(\mbf x^{k-1})$\\
\bottomrule
\end{tabular*}
\end{center}

When $p = q = 1$, Algorithm 2 is the same as Algorithm 1. By varying the row partition parameter $s$ and the column partition parameter $t$, we obtain the following special cases.  
\bit 
\item[(i)] $s=m$, $t=n$
$$\mbf x^k=\mbf x^{k-1}-\alpha_k\wh{\mbf A}_{i,j}\l(\frac{\wh{\mbf A}_{i,:}\mbf x^{k-1}}{p^2}-\frac{\wh{\mbf b}_i}{pq}-\frac{1-p}{p^2}\wh{\mbf A}_{i,j}(\mbf I_{:,j})^\rmt\mbf x^{k-1}\r)\mbf I_{:,j}.$$
\item[(ii)] $s=m$, $t=1$ 
$$\mbf x^k=\mbf x^{k-1}-\alpha_k\l((\wh{\mbf A}_{i,:})^\rmt\l(\frac{\wh{\mbf A}_{i,:}\mbf x^{k-1}}{p^2}- \frac{\wh{\mbf b}_i}{pq}\r)-\frac{1-p}{p^2}\diag\l((\wh{\mbf A}_{i,:})^\rmt\wh{\mbf A}_{i,:}\r)\mbf x^{k-1}\r).$$
\item[(iii)] $s=1$, $t=n$
$$\mbf x^k=\mbf x^{k-1}-\alpha_k(\wh{\mbf A}_{:,j})^\rmt\l(\frac{\wh{\mbf A}\mbf x^{k-1}}{p^2}-\frac{\wh{\mbf b}}{pq}-\frac{1-p}{p^2}\wh{\mbf A}_{:,j}(\mbf I_{:,j})^\rmt\mbf x^{k-1}\r)\mbf I_{:,j}.$$
\item[(iv)] $s=1$, $t=1$
$$\mbf x^k=\mbf x^{k-1}-\alpha_k\l(\wh{\mbf A}^\rmt\l(\frac{\wh{\mbf A}\mbf x^{k-1}}{p^2}-\frac{\wh{\mbf b}}{pq}\r)-\frac{1-p}{p^2}\diag\l(\wh{\mbf A}^\rmt\wh{\mbf A}\r)\mbf x^{k-1}\r).$$
\eit
We remark that at each step the cases (i) and (ii) only use the data $\wh{\mbf A}_{i,:}$ and $\wh{\mbf b}_i$ to update the iterate. This is particularly appropriate for those problems with  extremely large $m\times n$ matrix $\mbf A$ where it is not possible to load more than one row of $\mbf A$ due to memory constraints.

\subsection{Convergence analysis}

First, we will prove two useful properties of the update function $\mbf g(\mbf x)$ defined by (\ref{gx}). Lemma \ref{lem-g} shows that $\mbf g(\mbf x)$ is a stochastic gradient of the objective function $\dsp\frac{1}{2st}\|{\bf Ax-b}\|_2^2$. Lemma \ref{lem2} provides a uniformly upper bound of the expected norm of $\mbf g(\mbf x)$. 

\begin{lemma}\label{lem-g} For any fixed $\mbf x\in\mbbr^n$, let \beq\label{gx}\mbf g(\mbf x)=\mbf I_{:,\mcalj}(\wh{\mbf A}_{\mcali,\mcalj})^\rmt\l(\dsp\frac{\wh{\mbf A}_{\mcali,:}\mbf x}{p^2}-\frac{\wh{\mbf b}_\mcali}{pq}\r)-\dsp\frac{1-p}{p^2}\diag\l(\mbf I_{:,\mcalj}(\wh{\mbf A}_{\mcali,\mcalj})^\rmt\wh{\mbf A}_{\mcali,:}\r)\mbf x.\eeq We have $$\mbbe\bem \mbf g(\mbf x)\eem=\frac{1}{st}\mbf A^\rmt(\mbf A\mbf x-\mbf b).$$
\end{lemma}
\proof
Let $\mbbe_\delta\bem\cdot\eem$ denote the expectation with respect to the random binary mask (i.e., $\delta_{i,j}$ and $\delta_i$, in total $2^{m(n+1)}$ possibilities) and $\mbbe_{(\mcali,\mcalj)}\bem\cdot\eem$ denote the expectation with respect to the random selection ($st$ possibilities) of the pair of index sets. Then the full expected value $\mbbe\bem\cdot\eem$ satisfies $$\mbbe\bem\cdot\eem=\mbbe_{\delta}\mbbe_{(\mcali,\mcalj)}\bem\cdot\eem.$$
We recall that all $\delta_{i,j}$ and $\delta_i$ are independent. By straightforward calculations, we have 
\begin{align} \mbbe\bem \mbf I_{:,\mcalj}(\wh{\mbf A}_{\mcali,\mcalj})^\rmt\wh{\mbf A}_{\mcali,:} \eem & =\mbbe_{\delta}\mbbe_{(\mcali,\mcalj)}\bem \mbf I_{:,\mcalj}(\wh{\mbf A}_{\mcali,\mcalj})^\rmt\wh{\mbf A}_{\mcali,:}\eem \nn \\ &=\mbbe_{\delta}\mbbe_{(\mcali,\mcalj)}\bem \mbf I_{:,\mcalj}(\mbf I_{:,\mcalj})^\rmt\wh{\mbf A}^\rmt\mbf I_{:,\mcali}(\mbf I_{:,\mcali})^\rmt\wh{\mbf A}\eem \nn \\ &= \frac{1}{st}\mbbe_\delta\bem\wh{\mbf A}^\rmt\wh{\mbf A}\eem \nn\\ &=\frac{1}{st}\mbbe_\delta\bem (\wh{\mbf A}_{:,1})^\rmt\wh{\mbf A}_{:,1} & (\wh{\mbf A}_{:,1})^\rmt\wh{\mbf A}_{:,2} &\cdots & (\wh{\mbf A}_{:,1})^\rmt\wh{\mbf A}_{:,n}\\ (\wh{\mbf A}_{:,2})^\rmt\wh{\mbf A}_{:,1} & (\wh{\mbf A}_{:,2})^\rmt\wh{\mbf A}_{:,2} &\cdots & (\wh{\mbf A}_{:,2})^\rmt\wh{\mbf A}_{:,n}\\ \vdots & \vdots & \ddots & \vdots\\ (\wh{\mbf A}_{:,n})^\rmt\wh{\mbf A}_{:,1} & (\wh{\mbf A}_{:,n})^\rmt\wh{\mbf A}_{:,2} &\cdots & (\wh{\mbf A}_{:,n})^\rmt\wh{\mbf A}_{:,n}\eem\nn\\ &=\frac{1}{st} \bem p({\mbf A}_{:,1})^\rmt {\mbf A}_{:,1} & p^2({\mbf A}_{:,1})^\rmt {\mbf A}_{:,2} &\cdots & p^2({\mbf A}_{:,1})^\rmt{\mbf A}_{:,n}\\ p^2({\mbf A}_{:,2})^\rmt{\mbf A}_{:,1} & p({\mbf A}_{:,2})^\rmt{\mbf A}_{:,2} &\cdots & p^2({\mbf A}_{:,2})^\rmt{\mbf A}_{:,n}\\ \vdots & \vdots & \ddots & \vdots\\ p^2({\mbf A}_{:,n})^\rmt{\mbf A}_{:,1} & p^2({\mbf A}_{:,n})^\rmt{\mbf A}_{:,2} &\cdots & p({\mbf A}_{:,n})^\rmt{\mbf A}_{:,n}\eem\nn \\ \label{ATA} &=\frac{p^2}{st}\mbf A^\rmt\mbf A+\frac{p-p^2}{st}\diag(\mbf A^\rmt\mbf A). 
\end{align}
Similarly, we have \begin{align}\mbbe\bem\mbf I_{:,\mcalj}(\wh{\mbf A}_{\mcali,\mcalj})^\rmt\wh{\mbf b}_\mcali \eem &=\mbbe_\delta\mbbe_{(\mcali,\mcalj)}\bem\mbf I_{:,\mcalj}(\wh{\mbf A}_{\mcali,\mcalj})^\rmt\wh{\mbf b}_\mcali \eem \nn \\ &=\mbbe_\delta\mbbe_{(\mcali,\mcalj)}\bem\mbf I_{:,\mcalj}(\mbf I_{:,\mcalj})^\rmt\wh{\mbf A}^\rmt\mbf I_{:,\mcali}(\mbf I_{:,\mcali})^\rmt\wh{\mbf b} \eem \nn \\ &=\frac{1}{st}\mbbe_\delta\bem \wh{\mbf A}^\rmt\wh{\mbf b}\eem\nn \\ &=\frac{1}{st}\mbbe_\delta\bem \wh{\mbf A}^\rmt\eem\mbbe_\delta\bem \wh{\mbf b}\eem \nn\\ \label{ATb}&=\frac{pq}{st}\mbf A^\rmt\mbf b.\end{align}
Using (\ref{ATA}), we have \begin{align} \label{dATA}\mbbe\bem \diag\l(\mbf I_{:,\mcalj}(\wh{\mbf A}_{\mcali,\mcalj})^\rmt\wh{\mbf A}_{\mcali,:}\r)\eem=\diag\l(\mbbe\bem \l(\mbf I_{:,\mcalj}(\wh{\mbf A}_{\mcali,\mcalj})^\rmt\wh{\mbf A}_{\mcali,:}\r)  \eem\r)=\frac{p}{st}\diag(\mbf A^\rmt\mbf A).\end{align} 
Combining  (\ref{ATA}), (\ref{ATb}), and (\ref{dATA}) yields  $$\mbbe\bem \mbf g(\mbf x)\eem=\frac{1}{st}\mbf A^\rmt(\mbf A\mbf x-\mbf b).$$
This completes the proof. \endproof

\begin{lemma} \label{lem2}  For any fixed $\mbf x\in\mbbr^n$, let $\mbf g(\mbf x)$ be given as in (\ref{gx}). We have
\begin{align*}\mbbe\bem\|\mbf g(\mbf x)\|_2^2\eem &\leq \frac{2}{stp^2}\|\mbf A\|_\rmf^2\|\mbf A\mbf x-\mbf b\|_2^2  +\frac{2(1-q)}{stp^2q}\|\mbf A\|_\rmf^2\|\mbf b\|_2^2\\ &\quad +\frac{2(1-p)}{stp^3}\|\mbf A\|_\rmf^2\mbf x^\rmt\diag(\mbf A^\rmt\mbf A)\mbf x\\ &\quad +\frac{2(1-p)^2}{stp^3} \l\|\diag(\mbf A^\rmt\mbf A)\mbf x\r\|_2^2.
\end{align*}
\end{lemma}
\proof By straightforward calculations, we have \begin{align} \|\mbf g(\mbf x)\|_2^2& =\l\|\mbf I_{:,\mcalj}(\wh{\mbf A}_{\mcali,\mcalj})^\rmt\l(\dsp\frac{\wh{\mbf A}_{\mcali,:}\mbf x}{p^2}-\frac{\wh{\mbf b}_\mcali}{pq}\r)-\dsp\frac{1-p}{p^2}\diag\l(\mbf I_{:,\mcalj}(\wh{\mbf A}_{\mcali,\mcalj})^\rmt\wh{\mbf A}_{\mcali,:}\r)\mbf x\r\|_2^2\nn\\ & \leq\l(\l\|\mbf I_{:,\mcalj}(\wh{\mbf A}_{\mcali,\mcalj})^\rmt\l(\dsp\frac{\wh{\mbf A}_{\mcali,:}\mbf x}{p^2}-\frac{\wh{\mbf b}_\mcali}{pq}\r)\r\|_2+\l\|\dsp\frac{1-p}{p^2}\diag\l(\mbf I_{:,\mcalj}(\wh{\mbf A}_{\mcali,\mcalj})^\rmt\wh{\mbf A}_{\mcali,:}\r)\mbf x\r\|_2\r)^2\nn\\ \label{ine1} &\leq 2\l\|\mbf I_{:,\mcalj}(\wh{\mbf A}_{\mcali,\mcalj})^\rmt\l(\dsp\frac{\wh{\mbf A}_{\mcali,:}\mbf x}{p^2}-\frac{\wh{\mbf b}_\mcali}{pq}\r)\r\|_2^2+\frac{2(1-p)^2}{p^4} \l\|\diag\l(\mbf I_{:,\mcalj}(\wh{\mbf A}_{\mcali,\mcalj})^\rmt\wh{\mbf A}_{\mcali,:}\r)\mbf x\r\|_2^2\end{align} 
and 
\begin{align}
\l\|\mbf I_{:,\mcalj}(\wh{\mbf A}_{\mcali,\mcalj})^\rmt\l(\dsp\frac{\wh{\mbf A}_{\mcali,:}\mbf x}{p^2}-\frac{\wh{\mbf b}_\mcali}{pq}\r)\r\|_2^2 & =\l\|\mbf I_{:,\mcalj}(\wh{\mbf A}_{\mcali,\mcalj})^\rmt(\mbf I_{:,\mcali})^\rmt\l(\dsp\frac{\wh{\mbf A}\mbf x}{p^2}-\frac{\wh{\mbf b}}{pq}\r)\r\|_2^2 \nn \\ &\leq \l\|\mbf I_{:,\mcalj}(\wh{\mbf A}_{\mcali,\mcalj})^\rmt(\mbf I_{:,\mcali})^\rmt\r\|_\rmf^2\l\|\dsp\frac{\wh{\mbf A}\mbf x}{p^2}-\frac{\wh{\mbf b}}{pq}\r\|_2^2\nn\\ &\leq \l\|\mbf A_{\mcali,\mcalj}\r\|_\rmf^2\l\|\dsp\frac{\wh{\mbf A}\mbf x}{p^2}-\frac{\wh{\mbf b}}{pq}\r\|_2^2. \label{ine3} \end{align}
Further calculations give the expectation 
\begin{align}
& \hspace{5.2mm}\mbbe\bem\l\|\mbf A_{\mcali,\mcalj}\r\|_\rmf^2\l\|\dsp\frac{\wh{\mbf A}\mbf x}{p^2}-\frac{\wh{\mbf b}}{pq}\r\|_2^2\eem \nn \\ &  =\mbbe_\delta\mbbe_{(\mcali,\mcalj)}\bem\l\|\mbf A_{\mcali,\mcalj}\r\|_\rmf^2\l\|\dsp\frac{\wh{\mbf A}\mbf x}{p^2}-\frac{\wh{\mbf b}}{pq}\r\|_2^2\eem =\frac{1}{st}\|\mbf A\|_\rmf^2\mbbe_\delta\bem\l\|\dsp\frac{\wh{\mbf A}\mbf x}{p^2}-\frac{\wh{\mbf b}}{pq}\r\|_2^2\eem\nn \\ &=\frac{1}{st}\|\mbf A\|_\rmf^2\l(\frac{1}{p^4}\mbf x^\rmt\mbbe_\delta\bem\wh{\mbf A}^\rmt\wh{\mbf A}\eem\mbf x-\frac{2}{p^3q}\mbf x^\rmt\mbbe_\delta\bem\wh{\mbf A}^\rmt\wh{\mbf b}\eem+\frac{1}{p^2q^2}\mbbe_\delta\bem\wh{\mbf b}^\rmt\wh{\mbf b}\eem\r)\nn\\ &=\frac{1}{st}\|\mbf A\|_\rmf^2\l(\frac{1}{p^2}\mbf x^\rmt\mbf A^\rmt\mbf A\mbf x+\frac{1-p}{p^3}\mbf x^\rmt\diag(\mbf A^\rmt\mbf A)\mbf x-\frac{2}{p^2}\mbf x^\rmt\mbf A^\rmt\mbf b+\frac{1}{p^2q}\mbf b^\rmt\mbf b\r)\nn \\ &=\frac{1}{stp^2}\|\mbf A\|_\rmf^2\|{\bf Ax-b}\|_2^2+\frac{1-q}{stp^2q}\|\mbf A\|_\rmf^2\|\mbf b\|_2^2+\frac{1-p}{stp^3}\|\mbf A\|_\rmf^2\mbf x^\rmt\diag(\mbf A^\rmt\mbf A)\mbf x.\label{eq1}
\end{align}
It follows from \begin{align}\nn\diag\l(\mbf I_{:,\mcalj}(\wh{\mbf A}_{\mcali,\mcalj})^\rmt\wh{\mbf A}_{\mcali,:}\r) &=\diag\l(\mbf I_{:,\mcalj}(\mbf I_{:,\mcalj})^\rmt\wh{\mbf A}^\rmt\mbf I_{:,\mcali}(\mbf I_{:,\mcali})^\rmt\wh{\mbf A}\r)\\ &=\sum_{j=1}^n\mbf I_{:,j}(\mbf I_{:,j})^\rmt \l(\mbf I_{:,\mcalj}(\mbf I_{:,\mcalj})^\rmt\wh{\mbf A}^\rmt\mbf I_{:,\mcali}(\mbf I_{:,\mcali})^\rmt\wh{\mbf A}\r)\mbf I_{:,j}(\mbf I_{:,j})^\rmt\nn\\ \nn &= \sum_{j\in\mcalj}\mbf I_{:,j}(\mbf I_{:,j})^\rmt\wh{\mbf A}^\rmt\mbf I_{:,\mcali}(\mbf I_{:,\mcali})^\rmt\wh{\mbf A}\mbf I_{:,j}(\mbf I_{:,j})^\rmt\nn\\&= \sum_{j\in\mcalj}\mbf I_{:,j}(\wh{\mbf A}_{\mcali,j})^\rmt\wh{\mbf A}_{\mcali,j}(\mbf I_{:,j})^\rmt\nn\\&= \mbf I_{:,\mcalj}\diag\l((\wh{\mbf A}_{\mcali,\mcalj})^\rmt\wh{\mbf A}_{\mcali,\mcalj}\r)(\mbf I_{:,\mcalj})^\rmt\nn\end{align} 
that \beqs\mbf 0 \preceq\diag\l(\mbf I_{:,\mcalj}(\wh{\mbf A}_{\mcali,\mcalj})^\rmt\wh{\mbf A}_{\mcali,:}\r)\preceq \diag\l(\wh{\mbf A}^\rmt\wh{\mbf A}\r)\preceq \diag\l(\mbf A^\rmt\mbf A\r).\eeqs
This yields
\begin{align}
\hspace{5.2mm}\l\|\diag\l(\mbf I_{:,\mcalj}(\wh{\mbf A}_{\mcali,\mcalj})^\rmt\wh{\mbf A}_{\mcali,:}\r)\mbf x\r\|_2^2 &= \mbf x^\rmt \l(\diag\l(\mbf I_{:,\mcalj}(\wh{\mbf A}_{\mcali,\mcalj})^\rmt\wh{\mbf A}_{\mcali,:}\r)\r)^2\mbf x\nn\\ &\leq \mbf x^\rmt \diag\l(\mbf A^\rmt\mbf A\r)\diag\l(\mbf I_{:,\mcalj}(\wh{\mbf A}_{\mcali,\mcalj})^\rmt\wh{\mbf A}_{\mcali,:}\r)\mbf x.\label{ddd}
\end{align}
Then by (\ref{dATA}) and (\ref{ddd}), we have \begin{align}
\mbbe\bem\l\|\diag\l(\mbf I_{:,\mcalj}(\wh{\mbf A}_{\mcali,\mcalj})^\rmt\wh{\mbf A}_{\mcali,:}\r)\mbf x\r\|_2^2\eem &\leq\mbbe\bem  \mbf x^\rmt \diag\l(\mbf A^\rmt\mbf A\r)\diag\l(\mbf I_{:,\mcalj}(\wh{\mbf A}_{\mcali,\mcalj})^\rmt\wh{\mbf A}_{\mcali,:}\r)\mbf x\eem \nn\\ &\leq \frac{p}{st}\mbf x^\rmt\l(\diag\l(\mbf A^\rmt\mbf A\r)\r)^2\mbf x \nn\\ &=\frac{p}{st}\l\|\diag\l(\mbf A^\rmt\mbf A\r)\mbf x\r\|_2^2.\label{four}
\end{align}
Combining (\ref{ine1}), (\ref{ine3}),  (\ref{eq1}), and (\ref{four}) yields
\begin{align*} \mbbe\bem\|\mbf g(\mbf x)\|_2^2 \eem &\leq2\mbbe\bem\l\|\mbf A_{\mcali,\mcalj}\r\|_\rmf^2\l\|\l(\dsp\frac{\wh{\mbf A}\mbf x}{p^2}-\frac{\wh{\mbf b}}{pq}\r)\r\|_2^2\eem \nn\\ & \quad +\frac{2(1-p)^2}{p^4} \mbbe\bem\l\|\diag\l(\mbf I_{:,\mcalj}(\wh{\mbf A}_{\mcali,\mcalj})^\rmt\wh{\mbf A}_{\mcali,:}\r)\mbf x\r\|_2^2\eem\\&\leq \frac{2}{stp^2}\|\mbf A\|_\rmf^2\|\mbf A\mbf x-\mbf b\|_2^2  +\frac{2(1-q)}{stp^2q}\|\mbf A\|_\rmf^2\|\mbf b\|_2^2\\ &\quad +\frac{2(1-p)}{stp^3}\|\mbf A\|_\rmf^2\mbf x^\rmt\diag(\mbf A^\rmt\mbf A)\mbf x\\ &\quad +\frac{2(1-p)^2}{stp^3} \l\|\diag(\mbf A^\rmt\mbf A)\mbf x\r\|_2^2.
\end{align*}
This completes the proof. \endproof

Next, we give the main result of this paper, which shows the convergence behavior of Algorithm 2 with a constant step size.

\begin{theorem}\label{main} Let $\mbf x^k$ denote the $k$th iterate of Algorithm 2 applied to the linear least squares problem (\ref{ls}) with partially observed data (\ref{pd}). For a constant step size $\dsp0<\alpha<\frac{\sigma^2_{\rm min}(\mbf A)}{st\rho}$ (i.e., all $\alpha_k=\alpha$), it holds \begin{align*} \mbbe\bem\|\mbf x^k-\mbf A^\dag\mbf b\|_2^2\eem  &\leq \l(1-\frac{2\alpha\sigma_{\rm min}^2(\mbf A)}{st}+2\alpha^2\rho\r)^k\|\mbf x^0-\mbf A^\dag\mbf b\|_2^2+\frac{\alpha  C}{\sigma_{\rm min}^2(\mbf A)-\alpha st \rho},\end{align*} 
where $$\rho=\l\|\mbbe\bem\mbf B^\rmt \mbf B\eem\r\|_2, \quad \mbf B=\dsp\frac{1}{p^2} \mbf I_{:,\mcalj}(\wh{\mbf A}_{\mcali,\mcalj})^\rmt\wh{\mbf A}_{\mcali,:}-\frac{1-p}{p^2}\diag\l(\mbf I_{:,\mcalj}(\wh{\mbf A}_{\mcali,\mcalj})^\rmt\wh{\mbf A}_{\mcali,:}\r),$$ and 
\begin{align*}C & =\frac{2}{p^2}\|\mbf A\|_\rmf^2\|\mbf A\mbf A^\dag\mbf b-\mbf b\|_2^2  +\frac{2(1-q)}{p^2q}\|\mbf A\|_\rmf^2\|\mbf b\|_2^2\\ &\quad +\frac{2(1-p)}{p^3}\|\mbf A\|_\rmf^2(\mbf A^\dag\mbf b)^\rmt\diag(\mbf A^\rmt\mbf A)\mbf A^\dag\mbf b\\ &\quad +\frac{2(1-p)^2}{p^3} \l\|\diag(\mbf A^\rmt\mbf A)\mbf A^\dag\mbf b\r\|_2^2.\end{align*}
\end{theorem}

\proof By Lemma \ref{lem-g}, we have \beq\label{gd}\mbbe\bem \mbf g(\mbf x^{k-1})\ |\mbf x^{k-1}\eem=\frac{1}{st}\mbf A^\rmt (\mbf A\mbf x^{k-1}-\mbf b).\eeq By Lemma \ref{lem2}, we have \beq\label{lesc}\mbbe\bem\|\mbf g(\mbf A^\dag\mbf b)\|_2^2\eem\leq \frac{C}{st}.\eeq
Straightforward calculations yield \begin{align} \|\mbf x^k- \mbf A^\dag\mbf b\|_2^2 & =\|\mbf x^{k-1}-\mbf A^\dag\mbf b-\alpha \mbf g(\mbf x^{k-1})\|_2^2 \nn\\ &= \|\mbf x^{k-1}-\mbf A^\dag \mbf b\|_2^2-2\alpha(\mbf x^{k-1}-\mbf A^\dag\mbf b)^\rmt \mbf g(\mbf x^{k-1})+\alpha^2\|\mbf g(\mbf x^{k-1})\|_2^2\nn\\ &\leq \|\mbf x^{k-1}-\mbf A^\dag\mbf b\|_2^2-2\alpha(\mbf x^{k-1}-\mbf A^\dag\mbf b)^\rmt \mbf g(\mbf x^{k-1})\nn\\  &\quad +\alpha^2(\|\mbf g(\mbf x^{k-1})-\mbf g(\mbf A^\dag\mbf b)\|_2+\|\mbf g(\mbf A^\dag\mbf b)\|_2)^2\nn\\ &\leq \|\mbf x^{k-1}-\mbf A^\dag\mbf b\|_2^2-2\alpha(\mbf x^{k-1}-\mbf A^\dag\mbf b)^\rmt \mbf g(\mbf x^{k-1})\nn\\  &\quad +2\alpha^2\|\mbf g(\mbf x^{k-1})-\mbf g(\mbf A^\dag\mbf b)\|_2^2+2\alpha^2\|\mbf g(\mbf A^\dag\mbf b)\|_2^2, \label{xab}\end{align}
and \begin{align}\mbbe\bem\|\mbf g(\mbf x^{k-1})-\mbf g(\mbf A^\dag\mbf b)\|_2^2\ |\mbf x^{k-1}\eem & =\mbbe\bem(\mbf x^{k-1}-\mbf A^\dag\mbf b)^\rmt\mbf B^\rmt\mbf B(\mbf x^{k-1}-\mbf A^\dag\mbf b)\ |\mbf x^{k-1}\eem\nn \\ &\leq \l\|\mbbe\bem\mbf B^\rmt\mbf B\eem\r\|_2\|\mbf x^{k-1}-\mbf A^\dag \mbf b\|_2^2\nn\\ &=\rho\|\mbf x^{k-1}-\mbf A^\dag \mbf b\|_2^2.\label{gdiff}\end{align}
Combining (\ref{gd}), (\ref{lesc}), (\ref{xab}), and (\ref{gdiff}) yields \begin{align*} \mbbe\bem \|\mbf x^k- \mbf A^\dag\mbf b\|_2^2 \ | \mbf x^{k-1} \eem & \leq \|\mbf x^{k-1}-\mbf A^\dag\mbf b\|_2^2-\frac{2\alpha}{st} (\mbf x^{k-1}-\mbf A^\dag\mbf b)^\rmt\mbf A^\rmt(\mbf A \mbf x^{k-1}-\mbf b)\\ &\quad +2\alpha^2\rho\|\mbf x^{k-1}-\mbf A^\dag\mbf b\|_2^2+2\alpha^2\frac{C}{st}\\& = \|\mbf x^{k-1}-\mbf A^\dag\mbf b\|_2^2-\frac{2\alpha}{st} (\mbf x^{k-1}-\mbf A^\dag\mbf b)^\rmt\mbf A^\rmt\mbf A( \mbf x^{k-1}-\mbf A^\dag\mbf b)\\ &\quad +2\alpha^2\rho\|\mbf x^{k-1}-\mbf A^\dag\mbf b\|_2^2+2\alpha^2\frac{C}{st}\\ & \leq \l(1-\frac{2\alpha\sigma_{\rm min}^2(\mbf A)}{st}+2\alpha^2\rho\r)\|\mbf x^{k-1}-\mbf A^\dag\mbf b\|_2^2+2\alpha^2\frac{C}{st}.\end{align*}
Therefore, by the law of total expectation, we have  \begin{align*} \mbbe\bem\|\mbf x^k-\mbf A^\dag\mbf b\|_2^2\eem & =\mbbe\bem\mbbe\bem\|\mbf x^k- \mbf A^\dag\mbf b\|_2^2 \ | \mbf x^{k-1}\eem\eem \\ &\leq   \l(1-\frac{2\alpha\sigma_{\rm min}^2(\mbf A)}{st}+2\alpha^2\rho\r)\mbbe\bem\|\mbf x^{k-1}-\mbf A^\dag\mbf b\|_2^2\eem+2\alpha^2\frac{C}{st}\\ & \leq \cdots\\ & \leq \l(1-\frac{2\alpha\sigma_{\rm min}^2(\mbf A)}{st}+2\alpha^2\rho\r)^k\|\mbf x^0-\mbf A^\dag\mbf b\|_2^2+\frac{2\alpha^2 \dsp\frac{C}{st}}{\dsp\frac{2\alpha\sigma_{\rm min}^2(\mbf A)}{st}-2\alpha^2\rho}\\ & = \l(1-\frac{2\alpha\sigma_{\rm min}^2(\mbf A)}{st}+2\alpha^2\rho\r)^k\|\mbf x^0-\mbf A^\dag\mbf b\|_2^2+\frac{\alpha C}{\sigma_{\rm min}^2(\mbf A)-\alpha st\rho}. 	
\end{align*}
This completes the proof. \endproof

When $p = q = 1$ and $\mbf b\in\ran(\mbf A)$, Theorem \ref{main} implies that $\mbf x^k$ in Algorithm 1 using sufficiently small positive constant $\alpha$ converges to $\mbf A^\dag\mbf b$.

\section{Numerical results}
In this section, we report numerical experiments to illustrate the theoretical results. In each experiment, all data are available. Partially observed data are realized by the mask independent random variables $\delta_{ij}$ and $\delta_i$. This makes the error $\|\mbf x^k-\mbf A^\dag\mbf b\|^2_2$ of Algorithm 2 computable. The initial guess $\mbf x^0=\mbf 0$ and the relative error $\|\mbf x^k-\mbf A^\dag\mbf b\|^2_2/\|\mbf A^\dag \mbf b\|_2^2$ is averaged over 10 trials.  All experiments are performed using MATLAB on a laptop with 2.7-GHz Intel Core i7 processor, 16 GB memory, and Mac operating system.

In Algorithm 2, for simplicity, we use the row partition $\{\mcali_i\}_{i=1}^s$ with $\dsp s=\lceil\frac{m}{\ell}\rceil$: \beqas\mcali_i&=&\{(i-1)\ell+1,(i-1)\ell+2,\ldots,i\ell\},\quad i=1,2,\ldots,s-1,\\ \mcali_s&=&\{(s-1)\ell+1,(s-1)\ell+2,\ldots,m\},\eeqas and the column partition $\{\mcalj_j\}_{j=1}^t$ with $\dsp t=\lceil\frac{n}{\tau}\rceil$: \beqas\mcalj_j&=&\{(j-1)\tau+1,(j-1)\tau+2,\ldots,j\tau\},\quad j=1,2,\ldots,t-1,\\ \mcalj_t&=&\{(t-1)\tau+1,(t-1)\tau+2,\ldots,n\}.\eeqas

In each experiment, the matrix $\mbf A$ is generated from a standard normal distribution: $$\tt A=randn(m,n),$$ so $\mbf A$ is a full column rank matrix with probability one. For the case $\mbf b\in\ran(\mbf A)$, we use $$\tt b=A*randn(n,1),$$ and for the case $\mbf b\notin\ran(\mbf A)$, we use  $$\tt b=A*randn(n,1)+null(A\mbox{'})*ones(m-n,1).$$

\begin{figure}[htb]
\centerline{\epsfig{figure=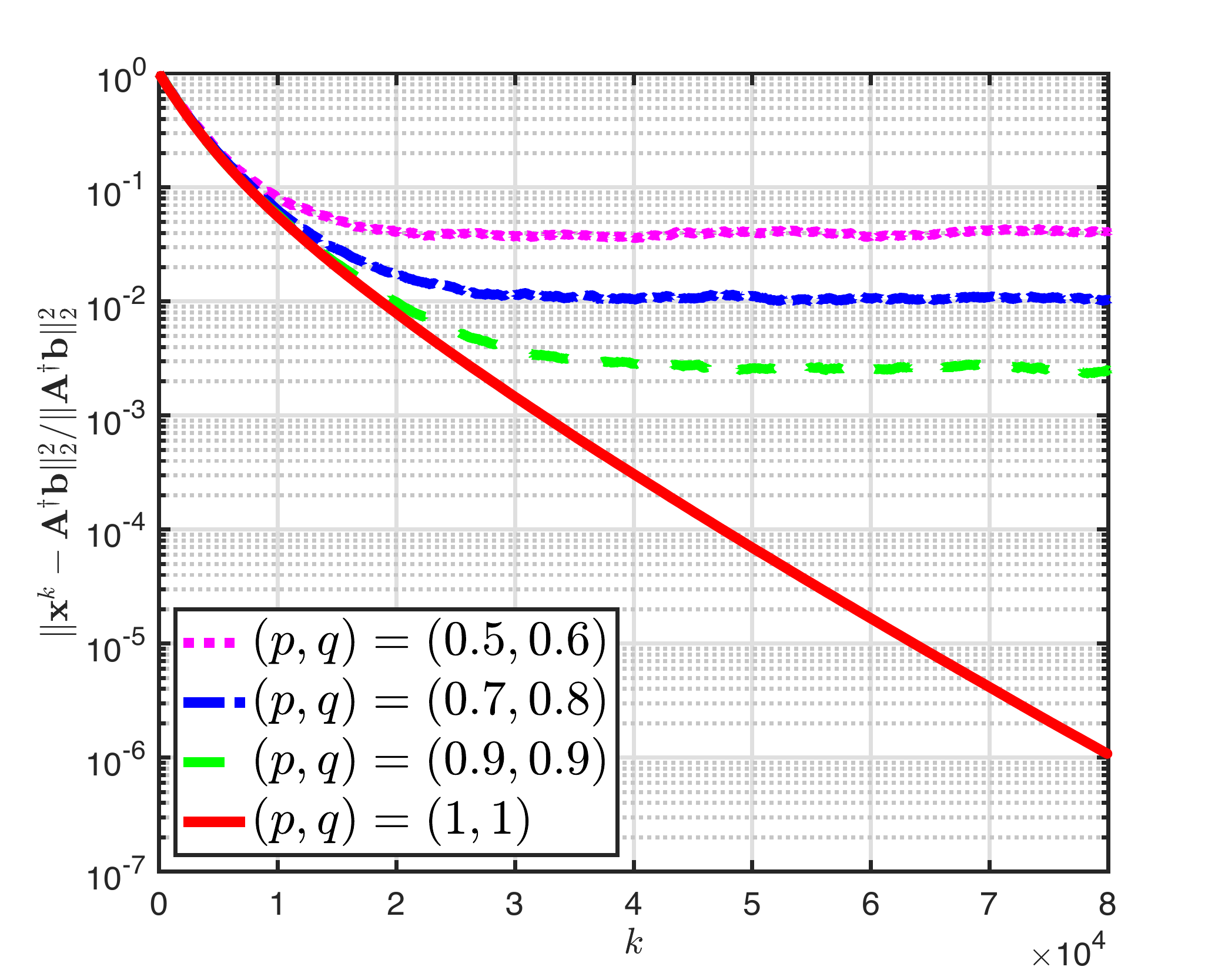,height=2.6in}\epsfig{figure=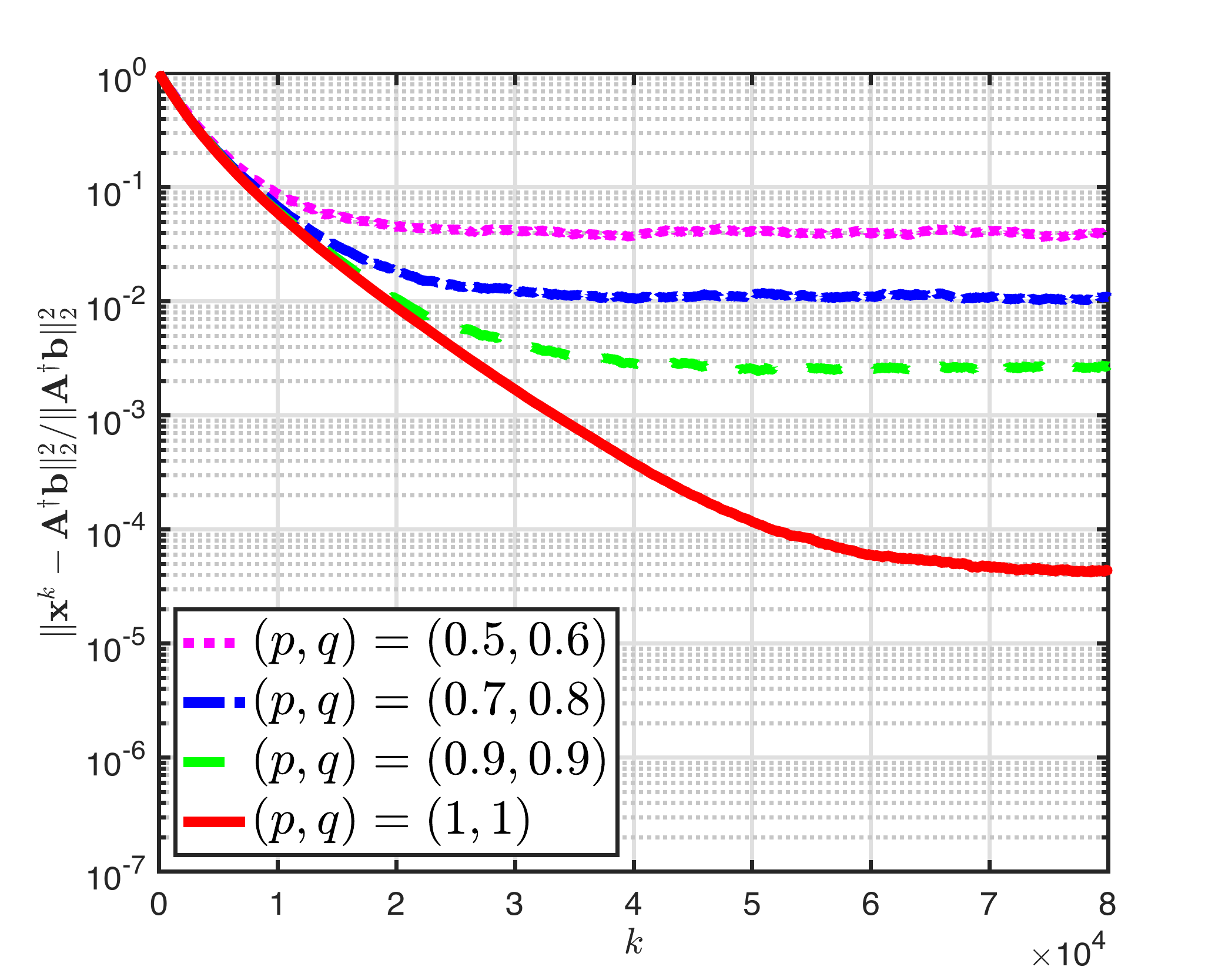,height=2.6in}}
\caption{The performance of Algorithm 2 using $\ell=2$, $\tau=n$, a constant step size $\alpha=10^{-4}$ and varied $p$ and $q$. Here, $m=1000$ and $n=200$. Left: $\mbf b \in\ran(\mbf A)$. Right: $\mbf b \notin\ran(\mbf A)$.}
\label{fig1}
\end{figure}

\begin{figure}[htb]
\centerline{\epsfig{figure=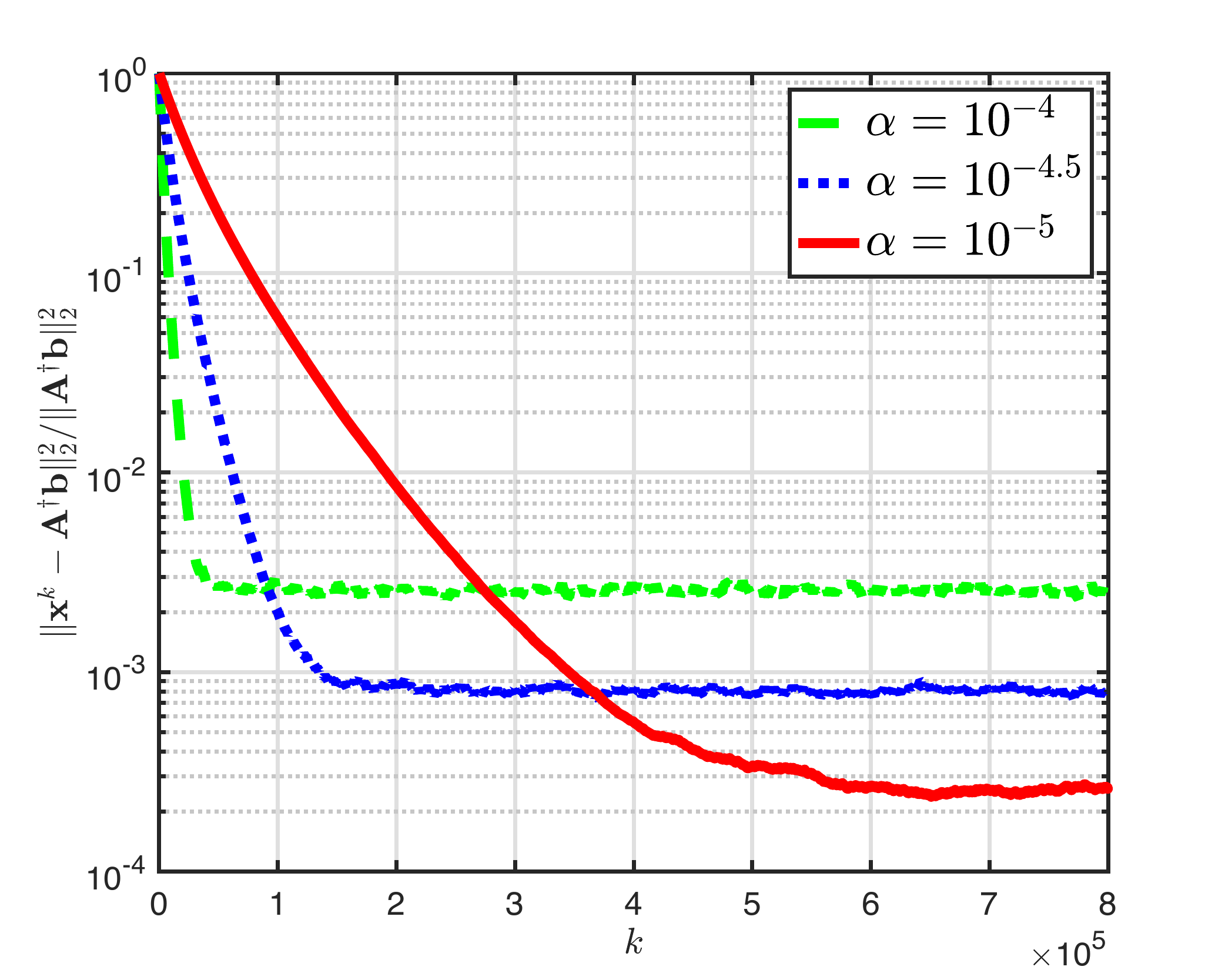,height=2.65in}\epsfig{figure=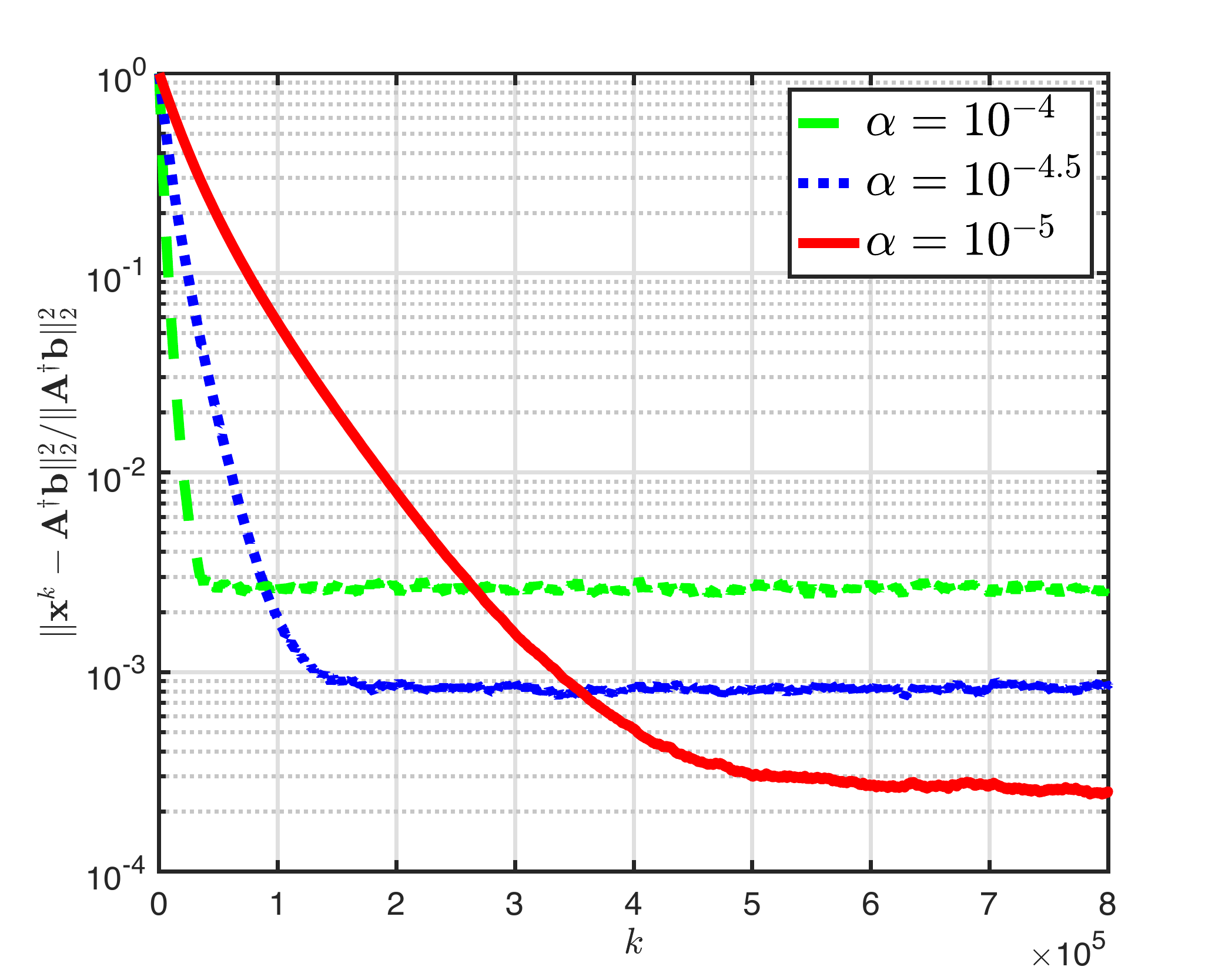,height=2.65in}}
\caption{The performance of Algorithm 2 using $\ell=2$, $\tau=n$, $p=0.9$, $q=0.9$, and different constant $\alpha$. Here, $m=1000$ and $n=200$. Left: $\mbf b \in\ran(\mbf A)$. Right: $\mbf b \notin\ran(\mbf A)$.}
\label{fig2}
\end{figure}

Figure \ref{fig1} shows the results of Algorithm 2 using $\ell=2$, $\tau=n$, a constant step size $\alpha=10^{-4}$ and varied proportions (i.e., $p$ and $q$) of available data. Figure \ref{fig2} shows the performance of Algorithm 2 using $\ell=2$, $\tau=n$, $p=0.9$, $q=0.9$, and different constant $\alpha$. These experimental results support the theoretical findings presented in Theorem \ref{main}. Using a constant step size, Algorithm 2 converges to some radius (proportional to $\alpha$) around the  solution. The proportions (i.e., $p$ and $q$) of available data affect the convergence horizon. In particular, as $p$ and $q$ decrease the size of the convergence horizon increases. When $p=q=1$ and $\mbf b\in\ran(\mbf A)$, Algorithm 2 behaves as DSBGS \cite{du2019doubl} does on the consistent linear system $\bf Ax = b$.


\begin{figure}[htb]
\centerline{\epsfig{figure=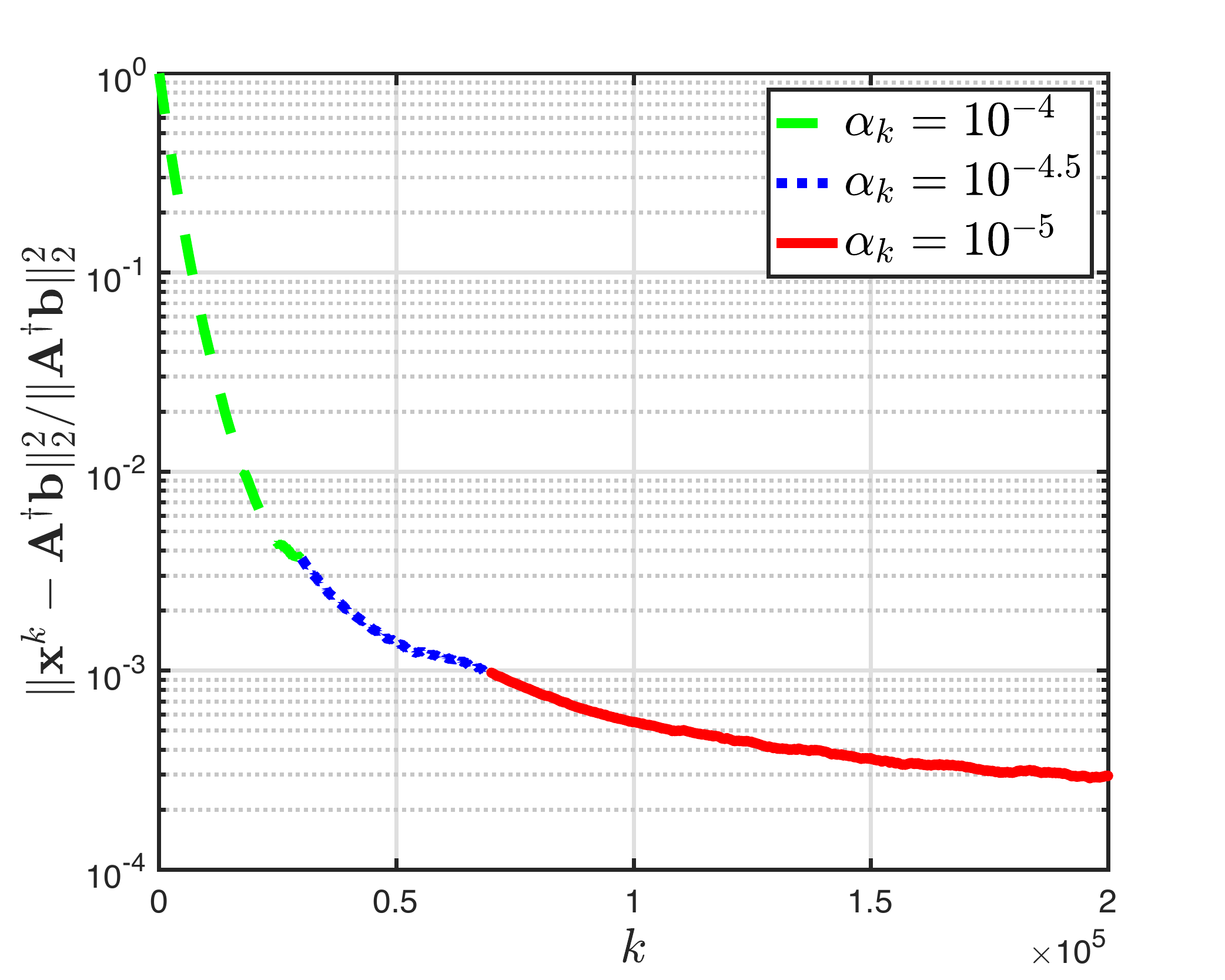,height=2.65in}\epsfig{figure=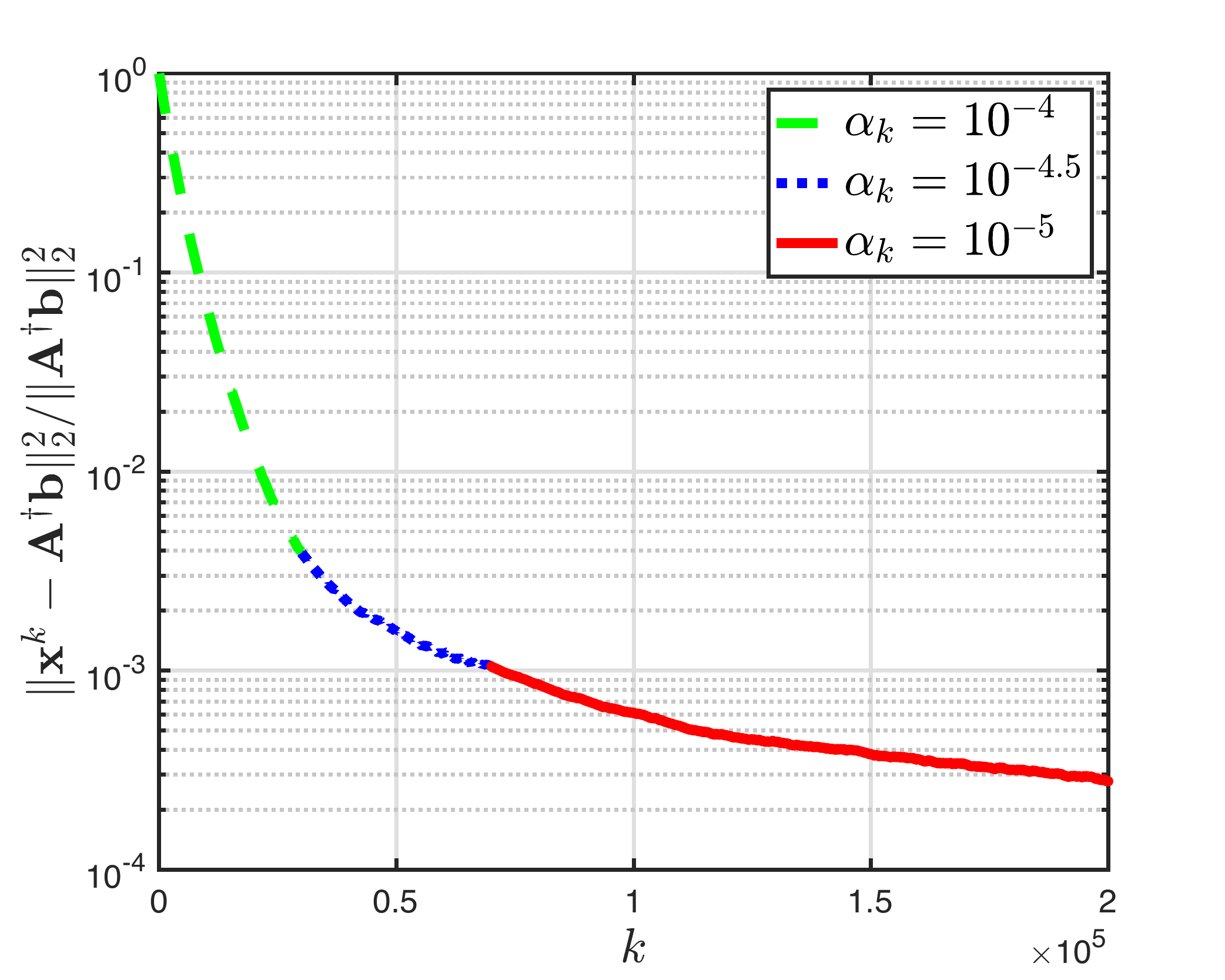,height=2.65in}}
\caption{The performance of Algorithm 2 using $\ell=2$, $\tau=n$, $p=0.9$, $q=0.9$, and updating $\alpha_k$. Here, $m=1000$ and $n=200$. Left: $\mbf b \in\ran(\mbf A)$. Right: $\mbf b \notin\ran(\mbf A)$.}
\label{fig3}
\end{figure}

Based on these numerical experiments, we can design a step size updating strategy: (i). choose pairs $\{(\beta_i,T_i)\}_{i=1}^K$ satisfying $\beta_1>\beta_2>\cdots>\beta_K>0$ and $T_1\leq T_2\leq \cdots\leq T_K$; (ii) use step size $\beta_1$ in the first $T_1$ iterations, and use step size $\beta_2$ in the following $T_2$ iterations, and so on. The performance of Algorithm  2 using this step size updating strategy with $\beta_1=10^{-4}$, $\beta_2=10^{-4.5}$, $\beta_3=10^{-5}$ and $T_1=3\times 10^4$, $T_2=4\times 10^4$, $T_3=1.3\times 10^5$, for the same data used in Figure \ref{fig2} is given in Figure \ref{fig3}. Compared with the constant step size strategy, the new strategy significantly reduces the number of iterations. 

\section{Concluding remarks}
We have proposed a stochastic gradient descent method for solving linear least squares problems with partially observed data. We prove that this method generates a sequence converging to some radius around the least squares solution. Numerical experiments on synthetic data illustrate the theoretical results. Finding appropriate step size selection strategies such as that used for Figure \ref{fig3}, and applying the resulting method on real world data should be valuable topics in the future study.

\section*{Acknowledgments}
The research of the first author was supported by the National Natural Science Foundation of China (No.11771364) and the Fundamental Research Funds for the Central Universities (No.20720180008).

\end{document}